\documentclass[12pt,a4paper]{article}
\usepackage{amsmath}
\usepackage{amssymb}
\usepackage{graphicx} 
\usepackage[sort,square,comma, numbers]{natbib}
\usepackage{hyperref}
\usepackage{anysize}
\marginsize{1.5cm}{1.5cm}{1cm}{1cm}

\newtheorem{proposition}{Proposition}[section]
\newtheorem{theorem}{Theorem}[section]
\newtheorem{corollary}{Corollary}[section]
\newtheorem{problem}{Problem}

\newtheorem{conjecture}{Conjecture}[section]
\newtheorem{observation}{Observation}[section]

\newenvironment{proof}[1][{}]%
{\par\noindent\textbf{Proof#1.}\quad}%
{\unskip\nobreak\hfil\penalty50\hskip2em\hbox{}%
\nobreak$\Box$\parfillskip=0pt\finalhyphendemerits=0\par}

\newcommand{\e}{\epsilon}

\newcommand{\Q}{\mathbb{Q}}
\newcommand{\R}{\mathbb{R}}
\newcommand{\Z}{\mathbb{Z}}
\newcommand{\N}{\mathbb{N}}
\newcommand{\Sp}{\mathrm{Spec}}
\newcommand{\cs}{\mathrm{cs}}
\newcommand{\charac}{\Psi}
\newcommand{\gold}{\varphi}
\newcommand{\EndProof}{\hspace{\stretch{1}} $\Box$}

\title{On graphs with cyclic defect or excess}

\author{
Charles Delorme\thanks{\texttt{cd@lri.fr}}\\
\emph{\small Laboratoire de Recherche en Informatique}\\
\emph{\small Universit\'e Paris-Sud}
\and
Guillermo Pineda-Villavicencio\thanks{\texttt{work@guillermo.com.au}}\\
\emph{\small Centre for Informatics and Applied Optimization}\\
\emph{\small University of Ballarat}\\
}

\begin{document}
\maketitle
\begin{abstract}

The Moore bound constitutes both an upper bound on the order of a graph
of maximum degree $d$ and diameter $D=k$ and a lower bound
on the order
of a graph of minimum degree $d$ and odd girth $g=2k+1$.

Graphs missing or exceeding the Moore bound by $\e$ are called {\it graphs with defect or excess $\e$}, respectively.

 Regular graphs with defect $\e$ satisfy the equation $G_{d,D}(A) = J_n + B$,
and regular graphs with excess $\e$ satisfy the equation $G_{d,\lfloor g/2 \rfloor}(A) = J_n - B$,
where $A$ denotes the adjacency matrix of the graph in question, $n$ its order,
$J_n$ the $n\times n$ matrix whose entries are all 1's,
$B$ a matrix with the row and column
sums equal to $\e$,
and  $G_{d,k}(x)$ a polynomial with integer coefficients such that
the matrix $G_{d,k}(A)$ gives the number of paths
of length at most $k$ joining each pair of vertices in the graph.

For {\it Moore graphs} (graphs with $\e=0$), the matrix $B$ is the null matrix. For graphs with defect or excess 1, $B$ can be considered as the adjacency matrix of a matching with $n$ vertices, while for graphs with defect or excess 2,  $B$ can be assumed to be the adjacency matrix of a union of vertex-disjoint cycles.

Graphs with defect 1 do not exist for any degree $\ge3$ and diameter $\ge2$, while graphs with excess 1 do not exist for any degree $\ge3$ and odd girth $\ge5$. However, graphs with defect or excess 2 represent a wide unexplored area.

Graphs with defect or excess 2 having the adjacency matrix of a cycle of order $n$ as the matrix $B$ are called \emph{graphs with cyclic defect or excess}; these graphs are the subject of our attention in this paper.

We obtain the following results about graphs with cyclic defect or excess. We prove the non-existence of
infinitely many such graphs, for example, graphs of any degree $\ge3$, diameter 3 or 4 and cyclic defect; and graphs of degree $\equiv0,2\pmod{3}$, girth 7 and cyclic excess. As the highlight of the paper we provide the asymptotic upper bound of $O(\frac{64}3d^{3/2})$ for the number of graphs of odd degree $d\ge3$ and cyclic defect or excess. This bound is in fact quite generous, and as a way of illustration, we also show that there are no graphs of degree 3 or 7, for diameter $\ge3$ and cyclic defect or for odd girth $\ge5$ and cyclic excess, nor any graphs of odd degree $\ge3$, girth $5$ or 9 and cyclic excess.

Actually, we conjecture
that, apart from the M\"obius ladder on 8 vertices,  no non-trivial graph of any degree $\ge 3$ and cyclic defect or excess exists.

To obtain our results we rely on algebraic methods, for instance, on the connection between the polynomial $G_{d,k}(x)$
and the classical Chebyshev polynomials, on eigenvalues techniques and on elements of algebraic number theory.

\end{abstract}

\noindent \textbf{Keywords:} Moore bound, Moore graph, defect, excess,
Chebyshev polynomial of the second kind, cyclic defect, cyclic excess, Pell equation.

\noindent\textbf{AMS  Subject Classification:} 05C12, 05C35, 05C50, 05C75.

\section{Introduction}

The terminology and notation used in this paper are standard and
consistent with that used in \cite{Die05}. Therefore, in this section we only settle the notation and terminology that could vary among texts.

The vertex set $V$ of a graph $\Gamma$ is denoted by $V(\Gamma)$, its
edge set by $E(\Gamma)$, its girth by $g(\Gamma)$, its adjacency matrix by
$A(\Gamma)$ and its diameter by $D(\Gamma)$; when there is no place for confusion,
we drop the symbol $\Gamma$. We often use the letter $n$ to denote the order of $\Gamma$.

The identity matrix of order $n$ is denoted by $I_n$, while by $J_n$ we denote
the $n\times n$ matrix whose entries are all 1's.

For a matrix $A$ the set formed by its $r+1$ distinct eigenvalues $\lambda_i$  with respective multiplicities
$m_i$ is called the {\it spectrum} of $A$ and is denoted by $\{[\lambda_0]^{m_0},\ldots,[\lambda_r]^{m_r}\}$.
The characteristic polynomial $\prod_{i=0}^r (x-\lambda_i)^{m_i}$ of $A$ is denoted by
$\charac(A,x)$. For a graph $\Gamma$, we often write $\charac(\Gamma,x)$ rather than
$\charac(A(\Gamma),x)$. We denote the eigenspace of $A$ corresponding to the eigenvalue
$\lambda$ by $E_\lambda(A)$.


We call a cycle of order $n$ an
{\it $n$-cycle} and denote it by $C_n$. If a graph $\Gamma$ is a  union of $m$ vertex-disjoint cycles,
we consider the multiset of their $r+1$ distinct lengths $l_i$ and respective multiplicities $m_i$,
and write that the {\it cycle structure} of $\Gamma$ is
$\cs(\Gamma)=\{[l_0]^{m_0},[l_1]^{m_1}\ldots [l_r]^{m_r}\}$ with
$m=\sum_{i=0}^r m_i$ and $n=\sum_{i=0}^r m_il_i$.

The degree of a polynomial $P$ is denoted by $\deg(P)$.
As it is customary, we denote the real Chebyshev polynomial of the second kind
by $U_m(x)$ \cite[pp.~3-5]{MH03}.
Recall that the polynomial $U_m(x)$, defined on $[-1,1]$, satisfies the following recurrence equations.
\begin{equation}\label{defU}
\begin{cases}
U_0(x)=1\\
U_1(x)=2x\\
U_{m+2}(x)=2xU_{m+1}(x)-U_m(x) \text{~for~} m\ge0 \text{~and~} x\in [-1,1]
\end{cases}
\end{equation}

It is known that the Moore bound, denoted by $M_{d,k}$ and defined below,
represents both an upper bound on the order of a graph of maximum degree $d$ and
diameter $D= k$ and a lower bound on
the order of a  graph of minimum degree $d$ and odd girth $g= 2k+1$ \cite{Big93}.
\begin{eqnarray}\label{boundgraph}
M_{d,k} & = &  1+d+d(d-1)+\ldots+d(d-1)^{k-1}
\nonumber
\\ & = & \left\{ \begin{array}{ll} 1 + d \frac{(d -1)^{k}
- 1}{d -2}
 & \mbox{if $d > 2$} \\ 2k+1 & \mbox{if $d = 2$} \\
\end{array}\right.
\end{eqnarray}
Non-trivial \emph{Moore graphs}
(graphs whose order equals the Moore bound, with $k\ge2$ and $d\ge 3$)
exist only for $D=2$ (or equivalently, for $g=5$), in which case
$d=2$, 3, 7 and possibly 57 \cite{HS60,BI73}.

By virtue of the rarity of Moore graphs, it is important to consider graphs which are somehow close to the ideal Moore graphs. Graphs of maximum degree $d$, diameter $D=k$ and order
$M_{d,k} - \e$ are called \emph{$(d,D,-\epsilon)$-graphs}, where the
parameter $\e$ is called \emph{defect}. Graphs of minimum degree $d$,
odd girth $g=2k+1$ and order $M_{d,k} + \e$ are called \emph{$(d,g,+\e)$-graphs}, where the parameter $\e$ is called the \emph{excess}.

Graphs with defect or excess 1 were completely classified by Bannai and Ito
\cite{BI81}; for any degree $d\ge2$, the only graphs of defect 1 are the cycles on $2D$ vertices, while the only graphs of excess 1 are the {\it cocktail party graphs} (the complement of $d/2+1$ copies of $K_2$, with even $d$).

\begin{figure}
\begin{center}
\includegraphics[scale=1]{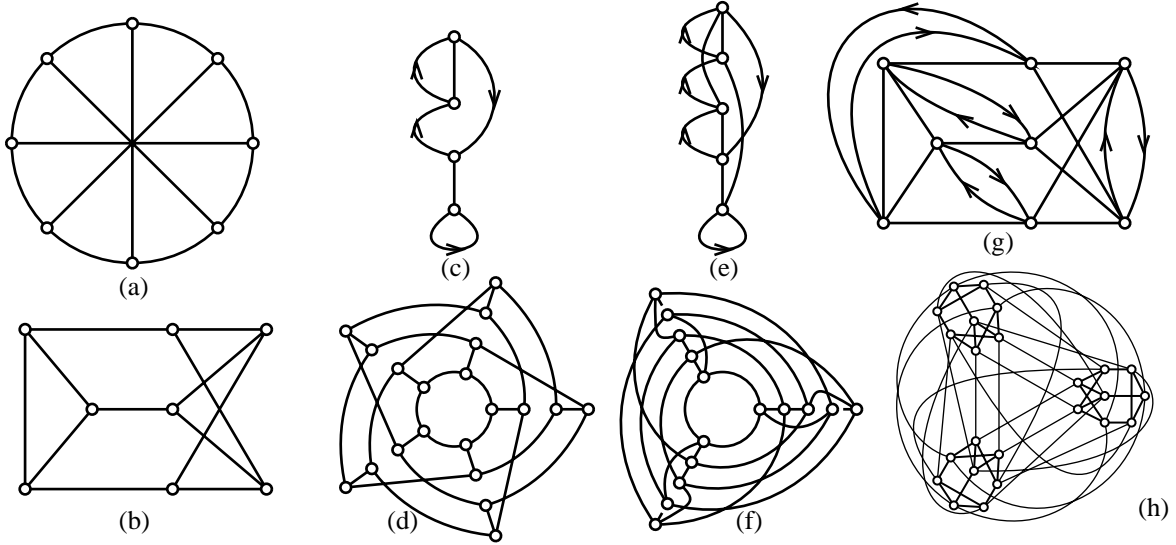}
\caption{All the non-trivial
known graphs with defect 2:
(a) the M\"obius ladder on 8 vertices, (b) the other $(3, 2,-2)$-graph, (c) a voltage graph of the unique $(3, 3,-2)$-graph, (d) the unique $(3, 3,-2)$-graph, (e) a voltage graph of the unique $(4, 2,-2)$-graph, (f) the unique $(4, 2,-2)$-graph, (g) a voltage graph of the unique $(5,2,-2)$ graph, and (h) the unique  $(5, 2,-2)$-graph.} \label{fig:Defect2}
\end{center}
\end{figure}

\begin{figure}
\begin{center}
\includegraphics[scale=1.1]{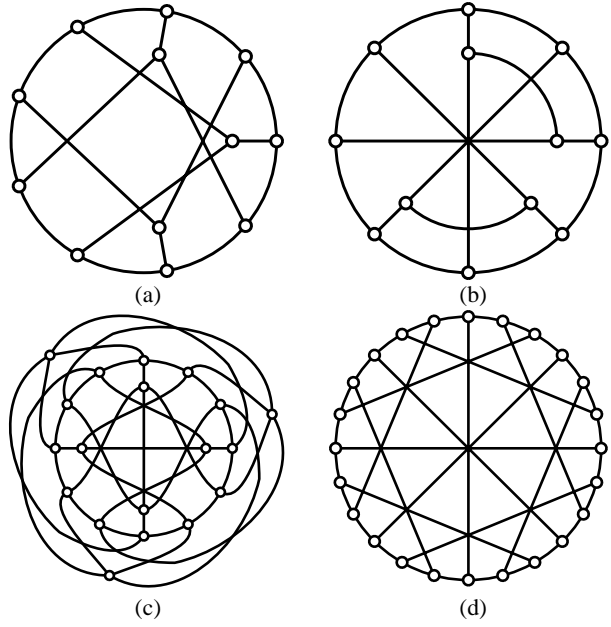}
\caption{All the non-trivial
known graphs with excess 2. (a) and (b) the only $(3,5,+2)$-graphs,  (c) the unique $(4, 5,+2)$-graph (the Robertson graph), and (d) the unique $(3, 7,+2)$-graph (the McGee graph).} \label{fig:Excess2}
\end{center}
\end{figure}

However, for $\e\ge 2$ the story is quite different. For maximum degree 2 and diameter $D\ge2$ the path of length $D$ is the only $(2,D,-2)$-graph. For degree $\ge3$ and diameter $D\ge2$  there are only 5 known graphs with defect 2, all of which are shown in Fig. \ref{fig:Defect2}. For degree 2 there is no graph with excess 2, while for degree $d\ge3$ and girth 3 the complement of the cycle $C_{d+3}$ is the only graph with excess 2. For degree $\ge3$ and odd girth $g\ge5$ there are only 4 graphs with excess 2 known at present (see Fig. \ref{fig:Excess2}).

For those familiar with the theory of voltage graphs (see \cite[Chapter 2]{GT87}), in Fig.~\ref{fig:Defect2} we present the $(3,3,-2)$-graph, the $(4,2,-2)$-graph and the $(5,2,-2)$-graph as lifts of voltage graphs. The $(3,2,-2)$-graph takes voltages on the group $\Z/5\Z$, while the $(4,2,-2)$-graph and the $(5,2,-2)$-graph take voltages on the group $\Z/3\Z$.  In all cases the undirected edges have voltage 0 and the directed edges have voltage 1.

It is worth mentioning that we gave an alternative voltage graph construction of a graph when this construction was simpler than the selected drawing of the graph. As principle failed for the $(3,2,-2)$-graphs, we omitted their respective voltage  graph representation.

It is not difficult to see that if $D=k\ge 2$ and $\e<1+(d-1)+\ldots+(d-1)^{k-1}$, a
$(d,D,-\e)$-graph must be $d$-regular. Similarly, if $g=2k+1\ge5$ and $\e<1+(d-1)+\ldots+(d-1)^{k-1}$, a
$(d,g,+\e)$-graph  must be $d$-regular.

Henceforth we consider graphs with defect or excess 2, and to avoid trivial cases, we only analyze graphs with degree $\ge 3$ and diameter $\ge 2$ for defect 2, and graphs with degree $\ge 3$ and girth $\ge 5$ for excess 2. Note that all these graphs must be regular.

In a graph $\Gamma$ with defect 2, if there are at least 2 paths of length
at most $D(\Gamma)$ from a vertex $v$ to a vertex $u$, then we say that $v$ is
a \emph{repeat} of $u$ (and \textit{vice versa}).
In this case we have  two repeats (not necessarily different)
for each vertex of $\Gamma$. Then, we define the {\it defect (multi)graph} of $\Gamma$ as the graph on $V(\Gamma)$, where two vertices are adjacent iff
one is a repeat of the other. Then, the defect graph is a union of vertex-disjoint cycles of length at least 2. Similarly, in a graph $\Gamma$ with excess 2, we define the \emph{excess graph}
of $\Gamma$ as the graph on $V(\Gamma)$,
where two vertices are adjacent iff they are at distance $D(\Gamma)$
(with $g(\Gamma)=2D(\Gamma)-1$).
Therefore, the excess graph is a union of vertex-disjoint cycles of length at least 3.

Next we present the cycle structure of the defect or excess graphs of the known non-trivial graphs with defect or excess 2.
\begin{description}
\item[Cyclic structure of graphs of defect 2] For the M\"obius ladder on 8 vertices $\cs=\{[8]^1\}$,  for the other $(3, 2,-2)$-graph $\cs=\{[3]^2,[2]^1\}$, for the unique $(3, 3,-2)$-graph $\cs=\{[5]^4\}$, for the unique $(4, 2,-2)$-graph $\cs=\{[6]^2,[3]^1\}$, and for the unique  $(5, 2,-2)$-graph $\cs=\{[3]^6,[2]^3\}$.
\item[Cyclic structure of graphs of excess 2] For the only $(3,5,+2)$-graphs (depicted in Fig.~\ref{fig:Excess2} as (a) and (b)) we have that (a) $\cs=\{[9]^1,[3]^1\}$ and (b) $\cs=\{[8]^1,[4]^1\}$, for the unique $(4, 5,+2)$-graph (the Robertson graph) $\cs=\{[3]^1,[12]^1,[4]^1\}$, and for the unique $(3, 7,+2)$-graph (the McGee graph) $\cs=\{[4]^6\}$.
\end{description}

For a graph $\Gamma$ of degree $d$ with adjacency matrix $A$, we define the
polynomials $G_{d,m}(x)$ for $x\in \R$:
\begin{equation}\label{defG}
\begin{cases}
G_{d,0}(x)=1\\
 G_{d,1}(x)=x+1\\
G_{d,m+1}(x)=xG_{d,m}(x)-(d-1)G_{d,m-1}(x) \text{~for~} m\ge 1
\end{cases}
\end{equation}

It is known that the entry $(G_{d,m}(A))_{\alpha ,\beta}$ counts the number of
paths of length at most $m$ joining the vertices $\alpha$ and $\beta$ in
$\Gamma$; see \cite{HS60,BI81,Si66}.

Regular graphs with defect $\e$ and order $n$ satisfy the equation
\begin{equation}
\label{eq:RepM} G_{d,D}(A) = J_n + B
\end{equation}
and regular graphs with excess $\e$ and order $n$ satisfy the equation
\begin{equation}
\label{eq:ExM}
G_{d,\lfloor g/2 \rfloor}(A) = J_n - B
\end{equation}
where $J_n$ is the $n\times n$ matrix whose entries are all 1's, and
$B$ is a matrix with the row and column sums equal to $\e$. The matrix $B$ is called the {\it defect} or {\it excess matrix} accordingly.

For Moore graphs, the matrix $B$ is the null matrix and $(x-d)G_{d,D}(x)$ is their minimal polynomial. For graphs with defect or excess 1, $B$ can be considered as the adjacency matrix of a matching with $n$ vertices \cite{BI81}. For a graph $\Gamma$ with defect or excess 2, the matrix $B$ is the adjacency matrix of the defect graph (respectively, of the excess graph). With a suitable labeling of $\Gamma$, $B$ becomes a direct sum of matrices representing cycles $C_l$ of length $l\ge2$ (respectively, $l\ge3$).

\[
A(C_2)=\begin{pmatrix} 0&2\\ 2&0\end{pmatrix}
\text{\hspace{10em}}A(C_l) =
\begin{pmatrix} 0&1&0&\dots&0&1\\
1&0&1&\dots&0&0\\
0&1&0&\dots&0&0\\
\hdotsfor[2]{6}\\
1&0&0&\dots&1&0
\end{pmatrix}
\]

\begin{figure}
\begin{center}
\includegraphics[scale=.9]{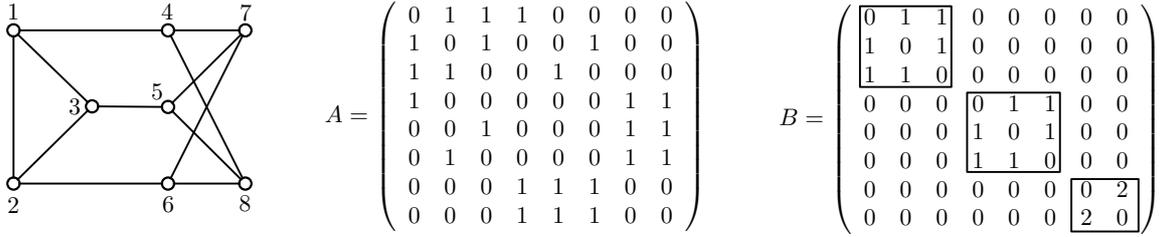}
\caption{Labelling of a $(3, 2,-2)$-graph that produces the desired structure of the corresponding defect matrix $B$.} \label{fig:322Graph}
\end{center}
\end{figure}

The previous point about the labelling of a graph $\Gamma$ is illustrated in Fig.~\ref{fig:322Graph}, where a $(3, 2,-2)$-graph is labelled such that the defect matrix $B$ displays the aforementioned structure.

For graphs with defect or excess 2, Equation~(\ref{eq:RepM}) has been studied for diameter $D=2$
\cite{CG08,NMP07,Fa87}, and
Equation~(\ref{eq:ExM}) has been studied for girths 5 and 7
\cite{Br67,Kov81,ES99}.

If $B$ is the adjacency matrix of a cycle of order $n$ (i.e.~$B=A(C_n)$), then the solution graphs of Equations~(\ref{eq:RepM}) and (\ref{eq:ExM}) are called \emph{graphs with cyclic defect} and \emph{graphs with cyclic excess}, respectively.

Among all
the known non-trivial graphs with defect or excess 2, only one has cyclic defect,
the M\"obius ladder on 8 vertices \cite{Fa87}, and none has cyclic excess.

In this paper we focus on graphs with cyclic defect or excess.
Basically, we deal with the following problems:

\begin{problem} \label{ProbRep} Classify the graphs of degree $d\ge3$, diameter $D\ge 2$
and order $n$ such that $G_{d,D}(A) =J_n +A(C_n)$.
\end{problem}
\begin{problem}\label{ProbEx} Classify the graphs of degree $d\ge3$, odd
girth $g\ge 5$ and order $n$ such that $G_{d,\lfloor g/2\rfloor}(A) =J_n - A(C_n)$.
\end{problem}
As Problem \ref{ProbRep} was completely settled for $D=2$ in \cite{Fa87}, from now on, we assume $D\ge3$.

The main result of the paper is the provision of the asymptotic upper bound of $O(\frac{64}3d^{3/2})$ for the number of graphs of odd degree $d\ge3$ and cyclic defect or excess. This bound turns out to be quite generous as our next results show. There are no graphs of degree 3 or 7, for diameter $\ge3$ and cyclic defect or for odd girth $\ge5$ and cyclic excess, nor any graphs of odd degree $\ge3$, girth $5$ or 9 and cyclic excess. Other non-existence outcomes are  the non-existence of graphs of any degree $\ge3$, diameter 3 or 4 and cyclic defect; and graphs of degree $\equiv0,2\pmod{3}$, girth 7 and cyclic excess.

To obtain our results we rely on algebraic methods, specifically on connections between the polynomials $G_{d,m}(x)$ and the classical Chebyshev polynomials of the second kind \cite{MH03}, on eigenvalue techniques, and on elements of algebraic number theory.

The rest of this paper is structured as follows.
In Section \ref{sec:Cycles} we provide some old and new combinatorial conditions for the existence of graphs with cyclic defect. In Section \ref{sec:AlgCond}
we present several algebraic approaches to analyze graphs with cyclic defect or excess, while Section \ref{sec:NonExRes} presents the main results of the paper. Finally, Section \ref{sec:Conclusion} summarizes our results and gives some concluding remarks.

\section{Combinatorial conditions for graphs with cyclic defect}\label{sec:Cycles}

Next we present some results about $(d,D,-2)$-graphs.

We denote by $\Theta_D$ the graph which is the union of three independent
paths of length $D$ with common endvertices.

\begin{proposition}{\cite[Lemma 2]{Jo92}}
Let $u$ be a vertex of a $(d,D,-2)$-graph $\Gamma$. Then
either:
\begin{enumerate}
\item[(i)] $u$ is a branch vertex of a $\Theta_D$ and
every cycle of length at most $2D$ in $\Gamma$ containing $u$ is
contained in this $\Theta_D$; or
\item[(ii)] $u$ is contained in one cycle of length $2D-1$ and no other
cycle of length at most $2D$; or
\item[(iii)] $u$ is
contained in exactly two cycles of length $2D$ and no other
cycle of length at most $2D$.
\end{enumerate}
\label{prop:Cycles}
\end{proposition}
\begin{corollary}
\label{cor:2DCycles}
Let $\Gamma$ be a $(d,D,-2)$-graph with cyclic defect. Then every vertex lies in exactly 2 cycles
of length $2D$.
\end{corollary}
\begin{corollary}\label{Ddivn}
The order $n$ of a  $(d,D,-2)$-graph with cyclic defect is a multiple of $D$.
\end{corollary}
\begin{proof}
By Corollary \ref{cor:2DCycles}, the number of $2D$-cycles in a  $(d,D,-2)$-graph with cyclic defect is $\frac{2n}{2D}$, and thus, the result follows.
\end{proof}
\begin{corollary}\label{congruences}
The allowed degrees for a $(d,D,-2)$-graph with cyclic defect are restricted
to some congruence classes modulo $D$.

When $D$ is even, $d$ is odd.

When $D$ is a power of an odd prime, $d-1$ is a multiple of $D$.

When $D\ge4$ is a power of 2,  $d-1$ is a multiple of $D/2$.
\end{corollary}
\begin{proof} If $2|D$, then $2|n$. As $n=M_{d,D}-2=-1+d(1+d-1+\ldots+(d-1)^{D-1})$, it follows that
$n\equiv d-1\pmod{d(d-1)}$, which implies $n\equiv d-1\pmod 2$, and thus, $2|(d-1)$.

Suppose that $D$ is a power of a prime $p$.

Suppose $d\equiv2\pmod p$. Then,
$1+(d-1)+\ldots+(d-1)^{D-1}\equiv D\pmod p$ and $n\equiv -1\pmod{p}$, which is incompatible
with $p|D|n$.

Suppose $d\not\equiv 2\pmod p$. By the little Fermat theorem \cite[p.~105]{Jac85I}, we have that $(d-1)^p\equiv d-1\pmod p$, and thus, that  $(d-1)^D\equiv d-1\pmod p$. Therefore, $n=-1+d\frac{(d-1)^D-1}{d-2}\equiv d-1\pmod p$. Also, since $n$ is a multiple of $D$, it follows that $d\equiv 1\pmod p$.

It remains to see what happens when $d\equiv1\pmod p$
for $D=p^r$ with $r>1$,
that is, $d=1+kp^s$ with $k\not\equiv 0\pmod p$.
 As $d=1+kp^s$, we have that $1+(d-1)+\ldots+(d-1)^{D-1}\equiv d\pmod{p^{s+1}}$. Therefore, it follows that $n=2kp^s\pmod{p^{s+1}}$.

Thus, to have $n\equiv 0\pmod D$ it is necessary that $s\ge r$ if $p$ is odd and $s\ge r-1$ if $p=2$. This completes the proof of the corollary. \end{proof}


\section{Algebraic conditions on the existence of graphs with cyclic defect or excess}
\label{sec:AlgCond}

We start this section by giving some known results.

If $B$ is the adjacency matrix of the $n$-cycle then its characteristic
polynomial $\charac(C_n,x)$
satisfies the following
\[
\charac(C_n,x)=\det(xI_n-B)=
\begin{cases}
(x-2)(x+2) (P_n(x))^2 & \text{~if~ $n$ is even}\\
(x-2) (P_n(x))^2 & \text{~if~ $n$ is odd}
\end{cases}
\]
where $P_n$ is a monic polynomial of degree $(n-2)/2$ if $n$ is even and  $(n-1)/2$ if $n$ is odd.

Recall that $x^n-1=\prod_{\ell|n}\Phi_\ell(x)$, where $\Phi_\ell(x)$ denotes the {\it $\ell$-th cyclotomic polynomial}\footnote{$\Phi_\ell(x)=\prod^{\phi(\ell)}_{m=1}(x-\xi_m)$, where $\{\xi_1,\xi_2,\ldots,\xi_{\phi(\ell)}\}$ denotes all $\ell$th primitive roots of unity, and $\phi(\ell)$ denotes the {\it Euler's totient function}, that is, the function giving  the number of positive integers  $\le\ell$ and relatively prime to $\ell$.}. The cyclotomic polynomial $\Phi_\ell(x)$  is an integer polynomial, irreducible over the field $\Q[x]$ of polynomials with rational coefficients, and {\it self-reciprocal} (that is, $x^{\phi(\ell)}\Phi_\ell(1/x)=\Phi_\ell(x)$). A consequence of $\Phi_\ell(x)$ being irreducible over $\Q[x]$ and self-reciprocal is that the degree of $\Phi_\ell(x)$ is even for $\ell\ge2$.

Using the previous facts on cyclotomic polynomials, we obtain the following factorization of $P_n(x)$: $P_n(x)=\prod_{3\le \ell|n}f_\ell(x)$, where  $f_\ell$ is an integer polynomial of degree $\phi(\ell)/2$ satisfying $x^{\phi(\ell)/2}f_\ell(x+1/x)=\Phi_\ell(x)$. Also, $f_\ell$ is irreducible over $\Q[x]$. In particular,  we have that $f_3(x)=x+1$, $f_4(x)=x$, $f_6(x)=x-1$,  $f_5(x)=x^2+x-1$, $f_8(x)=x^2-2$, $f_{12}(x)=x^2-3$,
$f_7(x)=x^3+x^2-2x-1$, $f_9(x)=x^3-3x+1$.

More concretely,
\begin{equation}\label{eq:BSpec}\Sp(B)=\begin{cases}
\{[2]^1,[2\cos{(\frac{2\pi}{n}\times 1)}]^{2},\ldots,
[2\cos{(\frac{2\pi}{n}\times \frac{n-2}{2})}]^{2},[-2]^1\}
&\text{if $n$ is even} \\
\{[2]^1,[2\cos{(\frac{2\pi}{n}\times 1)}]^{2},\ldots,
[2\cos{(\frac{2\pi}{n}\times \frac{n-1}{2})}]^{2}\}
&\text{if $n$ is odd}
\end{cases}
\end{equation}
It is also very well known that $\Sp(J_n)=\{[n]^1,[0]^{n-1}\}$.

Considering Equations \ref{eq:RepM} and \ref{eq:ExM}, we obtain that the eigenspace $E_n(J_n)$ equals both the eigenspace $E_d(A)$ and the eigenspace $E_2(B)$. Furthermore, for each eigenvalue $\lambda$ $(\ne d)$ of $A$, we have that $G_{d,D}(\lambda)$ is an eigenvalue $\mu$ $(\ne2)$
of $B$. In this case, we say that the eigenvalue $\lambda$ is paired with the eigenvalue $\mu$. Therefore, for each eigenvalue $\mu$ $(\ne 2)$ of $B$,
the eigenspace $E_\mu(B)$ contains the eigenspace of the eigenvalue of $A$ paired with $\mu$.

\begin{proposition}\label{prop:IntegerRootDef}
Let $A$ be the adjacency matrix of a $(d,D,-2)$-graph of order $n$.
If $n$ is even, then $A$ has a simple eigenvalue $\lambda$ such that
$\lambda$ is an integer root of the polynomial $G_{d,D}(x)+2$.
\end{proposition}

\begin{proof}  Consider Equations~(\ref{eq:RepM}) and
(\ref{eq:BSpec}). If $n$ is even, $-2$ is a simple eigenvalue of $B$,
and the eigenspace of $-2$ is spanned by the vector $u=(1,-1,1,-1,\ldots)^T$.
Let $\lambda$ be the simple eigenvalue of $A$ which is a root of $G_{d,D}(x)+2$.
Then, $u$ is also an eigenvector of $A$, implying that $\lambda$ must be integer.
\end{proof}

Let $\Gamma$ be a graph with cyclic defect.  If we substitute $y=G_{d,D}(x)$ into
$\charac(C_n,y)/(y-2)$,
we obtain a polynomial $F(x)$ of degree $(n-1)\times D$ such that $n-1$ of its roots are eigenvalues of $A$, and thus, $F(A)u=0$
for each eigenvector $u$ of $A$ orthogonal to the all-1 vector $\mathbf j$.

Setting $F_{\ell,d,D}(x):=f_\ell(G_{d,D}(x))$ we have
\[
F(x)=\begin{cases}
(G_{d,D}(x)+2)
\displaystyle{\prod_{\substack{\ell\; |\; n \\ \ell\ge 3}} (F_{\ell,d,D}(x))^2}
&\text{if $n$ is even}\\
\displaystyle{\prod_{\substack{\ell\; |\; n \\ \ell\ge 3}} (F_{\ell,d,D}(x))^2}
&\text{if $n$ is odd}.
\end{cases}
\]
\begin{observation}
\label{elemDiv1}
For each polynomial $f_\ell(x)$, where $\ell|n$ and $\ell\ge3$, the kernel of $f_\ell(B)$, denoted by $\ker(f_\ell(B))$, is formed by the direct sum of the eigenspaces associated with the roots of $f_\ell(x)$, and thus, $\ker(f_\ell(B))$ is a $\phi(\ell)$-dimensional space on $\Q[x]$.

Since $A$ commutes with $B$, we have that $\ker(f_\ell(B))$ is stable under the multiplication by $A$. Furthermore, as $B-G_{d,D}(A)$ is null on $\ker(f_\ell(B))$, it follows that $F_{\ell,d,D}(A)$ is null on $\ker(f_\ell(B))$ and that $\ker(F_{\ell,d,D}(A))$ is $\phi(\ell)$-dimensional on $\Q[x]$.

Consider a factor
$H(x)$ of $F_{\ell,d,D}(x)$. The kernel of $H(A)$ is stable under the multiplication by $B$, since $B-G_{d,D}(A)$ is null on $\ker(f_\ell(B))$. Thus, its dimension on $\Q[x]$ is either 0, or $\phi(\ell)/2$ or $\phi(\ell)$.

Hence, corresponding to the factor $f_\ell(x)$
of the minimal polynomial of $B$, the polynomial $F_{\ell,d,D}(x)$ has either 2 factors
of degree $\phi(\ell)/2$ or one factor of degree $\phi(\ell)$.
\end{observation}

By using Observation \ref{elemDiv1},
we obtain our first simple necessary
condition on the existence of graphs with cyclic defect.

\begin{proposition}\label{prop:irr}
For $D\ge3$ and $\ell\ge3$ such that $\ell|n$, if there is a $(d,D,-2)$-graph with cyclic
defect, then $F_{\ell,d,D}(x)$ must be reducible over $\Q[x]$.
\end{proposition}

\begin{proof}
Recall that $\deg(F_{\ell,d,D})=D\times \frac{\phi(\ell)}{2}$. If $F_{\ell,d,D}(x)$ is irreducible over $\Q[x]$, then all its roots must be eigenvalues of $A$. However, by Observation \ref{elemDiv1}, only $\phi(\ell)$ roots of $F_{\ell,d,D}(x)$ can be eigenvalues of $A$, a contradiction for $D\ge3$.
\end{proof}

Note that $\deg(F_{\ell,d,D})=D$ iff $\phi(\ell)=2$, and that $\phi(\ell)=2$ iff $\ell\in\{3,4,6\}$.
Thus, we have the following useful corollary.
\begin{corollary}\label{Cor:3,4,6}
Let $n$ be the order of a graph with cyclic defect and diameter $D\ge 3$. Then,
\begin{enumerate}
\item[(i)] if $n\equiv0\pmod3$ then $G_{d,D}(x)+1$ must be reducible over $\Q[x]$.
\item[(ii)] if $n\equiv0\pmod4$ then $G_{d,D}(x)$ must be reducible over $\Q[x]$.
\item[(iii)] if $n\equiv0\pmod6$ then $G_{d,D}(x)-1$ must be reducible over $\Q[x]$.
\end{enumerate}
\end{corollary}
\begin{proof}
 Knowing that $f_3(x)=x+1$, $f_4(x)=x$, and $f_6(x)=x-1$,
the result follows from Proposition \ref{prop:irr}.
\end{proof}

For $n\equiv0\pmod4$ we can even prove a result slightly stronger than the one of
Corollary \ref{Cor:3,4,6}.

Note that if $n\equiv0\pmod4$ then $d\equiv1\pmod2$.
As $n\equiv0\pmod4$, $0$ is an eigenvalue of $B$ with multiplicity 2.
 The vectors $u=(1,0,-1,0,1,\ldots)^T$
and $v=(0,1,0.-1,0,\ldots)^T$ form a basis of $E_{0}(B)$. As $A$ and $B$
commute, $Au\in E_{0}(B)$ and $Av\in E_{0}(B)$.
Therefore, we have that
\begin{equation}\label{eq:nEqui0(4)}
Au=\alpha u+\beta v \text{ and } Av=\delta u+\gamma v
\end{equation}
for some $\alpha,\beta,\delta,\gamma\in \Z$

Define a matrix $M$, called the {\it restriction of $A$ on $\ker(B)$}, as $\begin{pmatrix}\alpha&\delta\\ \beta&\gamma\end{pmatrix}$. Note that the characteristic polynomial $\Psi(M,x)$ of $M$ is the polynomial having as roots the two eigenvalues of $A$ paired with the eigenvalue 0 of $B$.

Let us consider  $u+v+\mathbf{j}$, where $\mathbf{j}$ is the all-1 vector. All components
of this sum are even. Thus, since all entries of $A$ are integers, $A(u+v+\mathbf{j})=(\alpha+\delta)u+(\beta+\gamma)v+d\mathbf{j}$ has
only even components. Consequently, $d+\alpha+\delta$ and $d+\beta+\gamma$ are even.

As $A$ is symmetric, $u ^ T A v=v ^T A u$ (recall that if $M_1$ and $M_2$ are matrices then
$(M_1M_2)^T=M_2^TM_1^T$).
Then, it follows that $u ^ T A v=u ^T (\delta u +\gamma v)=\frac n 2 \delta$ and that
 $v ^T A u = v^T (\alpha u + \beta v)= \frac n 2 \beta$, since
$u ^T u=\frac n 2 $,  $v^T v=\frac n 2 $ and $u ^T v =0$.
Thus, $\beta=\delta$ and $\alpha+\gamma\equiv 0\pmod{2}$.

In this way, we have obtained the following proposition.

\begin{proposition}\label{parity}
Let $A$ be the adjacency matrix of a graph with cyclic defect.
If  $n\equiv0\pmod{4}$
then the restriction of $A$ on the kernel of $B$ has an even trace. \EndProof
\end{proposition}

\begin{corollary}\label{cor:TrMDef} For $D=2$ the characteristic polynomial of the restriction of $A$ on the kernel of $B$ (i.e. $x^2+x+1-d$)
must be reducible over $\Q[x]$. \EndProof
\end{corollary}

The previous results on graphs with cyclic defect can be readily extended
to cover graphs with cyclic excess. Therefore, we limit ourselves to give
the results.

\begin{proposition}\label{prop:IntegerRootEx}
Let $A$ be the adjacency matrix of a $(d,g,+2)$-graph of order $n$.
If $n$ is even, there is a simple eigenvalue $\lambda$ of $A$ such that
$\lambda$ is an integer root of the polynomial $G_{d,\lfloor g/2\rfloor}(x) -2$. \EndProof
\end{proposition}

Let $\Gamma$ be a graph with cyclic excess.  Substituting
$y=-G_{d,\lfloor g/2\rfloor}(x)$ into
$\charac(C_n,y)/(y-2)$, we obtain a polynomial $F^*(x)$ of degree $(n-1)\times \lfloor g/2\rfloor$
such that $F^*(A)u=0$ for each vector $u$ orthogonal to the all-1 vector.
Setting $F^*_{\ell,d,\lfloor g/2\rfloor}(x):=f_\ell(-G_{d,\lfloor g/2\rfloor}(x))$ we have that
\[
F^*(x)=\begin{cases}
(-G_{d,\lfloor g/2\rfloor}(x)+2)
\displaystyle{\prod_{\substack{\ell\; |\; n \\ \ell\ge 3}}
(F^*_{\ell,d,\lfloor g/2\rfloor}(x))^2} &\text{if $n$ is even} \\
\displaystyle{\prod_{\substack{\ell\; |\; n \\ \ell\ge 3}}
(F^*_{\ell,d,\lfloor g/2\rfloor}(x))^2},&\text{if $n$ is odd}.
\end{cases}
\]
\begin{proposition}\label{prop:irrbis}
For $g\ge7$ and $\ell\ge3$ such that $\ell|n$, if there is a
$(d,g,+2)$-graph with cyclic excess and order $n$ then
$F^*_{\ell,d,\lfloor g/2\rfloor}(x)$ must be reducible over $\Q[x]$.
\end{proposition}

\begin{corollary}\label{Cor:3,4,6bis} Let $n$ be the order of a graph with cyclic excess. Then,
\begin{enumerate}
\item[(i)]
if $n\equiv0\pmod3$ then $G_{d,\lfloor g/2\rfloor}(x)-1$ must be reducible over $\Q[x]$.
\item[(ii)]
if $n\equiv0\pmod4$ then $G_{d,\lfloor g/2\rfloor}(x)$ must be reducible over $\Q[x]$.
\item[(iii)]
if $n\equiv0\pmod6$ then $G_{d,\lfloor g/2\rfloor}(x)+1$ must be reducible over
$\Q[x]$.
\end{enumerate}
\end{corollary}

\begin{proposition}\label{parityEx}
Let $A$ be the adjacency matrix of a graph with cyclic excess. If  $n\equiv0\pmod{4}$
then the restriction of $A$ on the kernel of $B$ has even trace. \EndProof
\end{proposition}

\begin{corollary}\label{cor:TrMEx} For $g=5$ and odd $d$, the characteristic polynomial of the restriction of $A$ on the kernel of $B$ (i.e. $x^2+x+1-d$)
must be reducible over $\Q[x]$. \EndProof
\end{corollary}

\subsection{Relations between the polynomials $G_{d,m}(x)$ and $U_{m}(x)$}
To establish some relations between the polynomials $G_{d,m}(x)$ and $U_m(x)$,
we make use of their respective generating functions (ordinary power series)
\[\displaystyle{P(x,t)=\frac{1+t}{1-xt+(d-1)t^2}}
\text{~and~}\displaystyle{Q(x,t)=\frac{1}{1-2xt+t^2}}.
\]
It is convenient to introduce $q:=\sqrt{d-1}>0$. Then, it follows that
$$
\begin{array}{c@{~=~}c@{~=~}c}
P(x,\frac{t}{q})&\sum_{m=0} q^{-m} G_{d,m}(x)t^m&\frac{1+t/q}{1-xt/q+t^2} \vspace{.2cm}\\
Q(\frac{x}{2q},t)&\sum_{m=0} U_{m}(\frac{x}{2q})t^m&\frac{1}{1-xt/q+t^2}
\end{array}
$$

Thus,
\begin{equation}\label{eq:GdD_UD}
G_{d,m}(x)=q^m U_{m}(\frac{x}{2q})+q^{m-1} U_{m-1}(\frac{x}{2q})
\end{equation}
Equation~(\ref{eq:GdD_UD}) allows us to establish some bounds
for the eigenvalues of a graph with cyclic defect (Proposition \ref{prop:EigBounds}).
\begin{proposition}\label{prop:EigBounds}
For a graph of degree $d\ge3$, diameter   $D\ge 2$ and cyclic defect, if $\beta$ is real and $|\beta|\le 2$,
then the roots of $G_{d,D}(x)+\beta$ are real and belong to the open interval $(-2\sqrt{d-1},2\sqrt{d-1})$.
\end{proposition}
\begin{proof} Set $q:=\sqrt{d-1}$, and  notice that $q^D\ge 2$, with equality only when $D=2$ and $d=3$.

From Equation (\ref{defU}) observe that $U_D(1)=D+1$ and $U_D(-1)=(-1)^D(D+1)$.
Then, using Equation (\ref{eq:GdD_UD}) we obtain that $$G_{d,D}(2q)=(D+1)q^D+Dq^{D-1}>2,$$ and that
$$G_{d,D}(-2q)=(-1)^D((D+1)q^D-Dq^{D-1}) \text{ has the sign of $(-1)^D$ and absolute value
$> 2.$}$$

We compute $G_{d,D}(2q\cos (t\pi/D))$ for $1\le t\le D-1$.
Since $U_{D-1}(\cos (t\pi/D))=0$, it follows that
$$G_{d,D}(2q\cos (t\pi/D))=q^DU_D(\cos (t\pi/D)) =q^D\; \frac{\sin(t(D+1)\pi/D)}{\sin(t\pi/D)}=(-1)^tq^D.$$

Hence, for any $|\beta|<2$, $d\ge3$ and $D\ge2$, with the exception of $\beta=2$, $D=2$ and $d=3$, each of the $D$ open intervals
$(2q\cos ((t+1)\pi/D),2q\cos (t\pi/D))$ with $0\le t\le D-1$ contains a root of $G_{d,D}(x)+\beta$ (by the Intermediate Value Theorem).

In the case $\beta=2$, $D=2$ and $d=3$,  the roots of $G_{3,2}(x)+2$ are 0 and $-1$,
which belong to ($-2\sqrt2$, $2\sqrt2$).
\end{proof}
By virtue of Proposition \ref{prop:EigBounds}, we can assume that every eigenvalue ($\ne d$) of $A$ has the form $2q\cos{\alpha}$ with $q:=\sqrt{d-1}$ and $0<\alpha<\pi$. In this case
\begin{equation}\label{eq:Gcos}
G_{d,m}(2q\cos{\alpha})=\frac{q\sin{(m+1)\alpha}+\sin{m\alpha}}{\sin\alpha}q^{m-1}, \textrm{ with $\sin\alpha\ne0$}
\end{equation}
As a corollary of Proposition \ref{prop:EigBounds}, we obtain a very useful necessary
condition on the
existence of graphs with cyclic defect and even order.

\begin{corollary}\label{cor:EvenN} If $n\equiv0\pmod{2}$ then
a graph with cyclic defect must have an integer eigenvalue $\lambda$ such that
$|\lambda|<2q$ and $G_{d,D}(\lambda)=-2$.
\end{corollary}
\begin{proof}  The corollary follows immediately from
Propositions \ref{prop:IntegerRootDef} and \ref{prop:EigBounds}.\end{proof}
Extending Proposition \ref{prop:EigBounds} and Corollary \ref{cor:EvenN}
to graphs with cyclic excess, we obtain the following assertions.

 \begin{proposition}\label{prop:EigBoundsEx}
For a graph of degree $d\ge3$, odd girth  $g\ge 5$ and cyclic excess, if $\beta$ is real and $|\beta|\le 2$,
then the roots of $G_{d,\lfloor g/2\rfloor}(x)+\beta$ are real
and belong to the open interval $(-2\sqrt{d-1},2\sqrt{d-1})$. \EndProof
\end{proposition}

\begin{corollary}\label{cor:EvenNEx} If $n\equiv0\pmod{2}$ then
a graph with cyclic excess must have an integer eigenvalue $\lambda$ such that
$|\lambda|<2q$ and $G_{d,\lfloor g/2\rfloor}(\lambda)=2$. \EndProof
\end{corollary}

\section{Results on graphs with cyclic defect or excess}
\label{sec:NonExRes}

\subsection{Graphs of diameter 4 and cyclic defect, and graphs of girth 9 and cyclic excess}

Here we basically prove the non-existence of graphs of degree $d\ge3$, diameter 4 and cyclic defect, or graphs of degree $d\ge3$, girth 9 and cyclic excess. Our proof is decomposed into three parts. We first show that the polynomial $G_{d,4}(x)$ has an integer root, then we find those values of $d$ making the existence of this root possible. Finally we show that for none of these values of $d$ the polynomial $G_{d,4}(x)\pm2$ has an integer root (contradicting Propositions \ref{prop:IntegerRootDef} and \ref{prop:IntegerRootEx}).

\begin{theorem}\label{theo:D4}
There is no regular graph of degree $d\ge3$, diameter 4 and cyclic defect, nor any regular graph of odd degree $d\ge3$, girth 9 and cyclic excess.
\end{theorem}

\begin{proof}
Considering graphs of diameter 4 and cyclic defect, from Corollary \ref{congruences} it follows that $d\equiv1\pmod2$ and that $n\equiv0\pmod4$, while for regular graphs of odd degree $d\ge3$ and cyclic excess it follows that  $n\equiv0\pmod4$. Set $a:=d-1$, then
$G_{a+1,4}=x^4+x^3-3ax^2-2ax+a^2$. By Corollary \ref{Cor:3,4,6} (for cyclic defect) and Corollary \ref{Cor:3,4,6bis} (for cyclic excess) the polynomial $G_{a+1,4}(x)$ must be reducible over $\Q[x]$, and thus, it must have a factor of degree at most 2. We first claim that for $a>1$ $G_{a+1,4}(x)$ must have an integer root.

{\bf Claim 1.} for $a>1$ $G_{a+1,4}(x)$ must have an integer root.

{\bf Proof of Claim 1.} We proceed by contradiction, assuming that there is a factorization of $G_{a+1,4}(x)$ into factors of degree 2 irreducible over $\Q[x]$. Then, from the roots $x_1,x_2,x_3,x_4$ of $G_{a+1,4}(x)$  we can obtain two sets, say $\{x_1,x_2\}$ and $\{x_3,x_4\}$, such that $x_1+x_2$, $x_1x_2$,  $x_3+x_4$ and $x_3x_4$ are all integers.

Using Vi\`ete's formulas we obtain that
\begin{align*}
\sigma _1&:=x_1+x_2+x_3+x_4=-1\\
\sigma _2&:=x_1x_2+x_1x_3+_1x_4+x_2x_3+x_2x_4+x_3x_4=-3a\\
\sigma _3&:=x_1x_2x_3+x_1x_2x_4+x_1x_3x_4+x_2x_3x_4=2a\\
\sigma _4&:=x_1x_2x_3x_4=a ^2
\end{align*}
 Therefore, we can compute the coefficients of the equation $p(y)=y^3-b_1y^2+b_2y-b_3$ with the 3 roots $y_1=x_1x_2+x_3x_4$, $y_2=x_1x_3+x_2x_4$ and $y_3=x_1x_4+x_2x_3$; indeed, we have that
 \begin{align*}b_1&=\sigma _2=-3a\\ b_2&=\sigma _1\sigma _3-4\sigma _4=-2a-4a ^2\\ b_3&=\sigma _2^2\sigma _4+\sigma _3^2-4\sigma _2\sigma _4=5a ^2 + 12 a ^3
 \end{align*}

Thus, we have to find integer solutions for $p(y)=y^3+3ay^2-4y a ^2  -12 a ^3  -2ay-5a ^2=0$. Discarding the uninteresting solution $y=a=0$, we may write $p(y)/a^2$ as
$(y-2a)u(u+1)-2u-1=0$, where $u:=2+y/a$ is rational and $y-2a\ne0$. This equation in $u$
has discriminant $(y-2a-2)^2+4(y-2a)=(y-2a)^2+4$, which can be a perfect square only if $y-2a=0$, a contradiction. Therefore, $p(y)$ cannot have integer roots, and the claim follows. \EndProof

Since $G_{a+1,4}(x)$ must have an integer root, we search the integer pairs $(x,a)$ such that $G_{a+1,4}(x)=0$.

The discriminant  $x^2(5x^2+8x+4)$
of the equation $G_{a+1,4}(x)=x^4+x^3-3ax^2-2ax+a ^2=0$ in $a$ is a perfect square iff
$5x^2+8x+4=t^2$; multiplying this equation by 5 and setting $z:=5x+4$, we obtain
\begin{equation}\label{eq:Pell}
z^2-5t^2=-4
\end{equation}
Equation (\ref{eq:Pell}) is closely related to the well-known Pell equation\footnote{While this equation is widely known as the Pell equation, there is no evidence that John Pell posed it. It seems that Euler was the causer of this confusion. See \cite[p.~4]{JW09} for more information.} (namely, $Z^2-PT^2=1$, where $Z,P,T\in \Z$).

The infinitely many  solutions $(z_m,t_m)$ of Equation (\ref{eq:Pell}) are given by $z_m=\pm L_{4m+3}$ and $t_m=\pm F_{4m+3}$, where $L_m$ and $F_m$ denote the $m$th Lucas number and $m$th Fibonacci number, respectively; see \cite[p.~64]{Sam70}. For all integers $m$  the recurrence equations of the Lucas and the Fibonacci numbers can be defined as follows.
\begin{equation}\label{eq:Luc&Fib}
\begin{cases}
L_0=2,L_1=1\\
L_{m+2}=L_{m+1}+L_m\\
\end{cases}\,\begin{cases}
F_0=0,F_1=1\\
F_{m+2}=F_{m+1}+F_m
\end{cases}
\end{equation}
If we set $\gold:=\frac{1+\sqrt 5}2$ (the so-called {\it golden ratio}), then $L_{4m+3}=\gold^{4m+3}-\gold^{-(4m+3)}$
and $F_{4m+3}=(\gold^{4m+3}+\gold^{-(4m+3)})/{\sqrt 5}$.
In order to retain integer values for $x$, we have
that $x_m=(-4+L_{4m+3})/5$, and thus, that $a_m=x_m(3x_m+2\pm t_m)/2$.

Set $r_m:=\gold^{4m+3}$, then $x_m=(r_m-1/r_m-4)/5$ and $t_m=(r_m+1/r_m)/\sqrt5$.

We first rule out the existence of graphs of diameter 4 and cyclic defect.

{\bf Claim 2.} There is no regular graph of degree $d\ge3$, diameter 4 and cyclic defect.

{\bf Proof of Claim 2.} For the aforementioned values of $a_m$, by Proposition \ref{prop:IntegerRootDef}, the polynomial  $G_{a_m+1,4}(x)+2$ must have an integer root. Our goal now is to prove that this is not the case.

From the two possible values for $a_m$ take $a_m=x_m(3x_m+2+ t_m)/2$.

Note that for any two integer values $u$ and $v$, $(G_{a_m+1,4}(u)-G_{a_m+1,4}(v))/(u-v)$ is an integer. Suppose that $u_m$ is an integer root of   $G_{a_m+1,4}(x)+2$, then, for $u_m$ and $x_m$ we have that $$\frac{G_{a_m+1,4}(u_m)-G_{a_m+1,4}(x_m)}{u_m-x_m}=\frac{-2}{u_m-x_m}$$ is an integer, which implies that $u_m-x_m=s_m=\pm 2$ or $\pm1$.

As a result, it follows that $$H(r_m):=r^3_m(G_{a_m+1,4}(u_m)+2)=r^3_m(G_{a_m+1,4}(x_m+s_m)+2)=0.$$ Note that $H(r_m)$ is a polynomial in $r_m$ of degree 6.

Investigating the real roots of $H(r_m)$ for each value of $s_m$, we see that their absolute values  lie between 0.05 and 9. But, since $r_m:=\gold^{4m+3}$, we have that for $m\ge1$, the values of $r_m$ are at least 29, and that for $m\le-3$, the values of $r_m$ lie between 0 and 0.01. For $m=-2$ it can be easily verified that $a_m=3$, which contradicts the fact that $a\equiv0\pmod2$. We have excluded the values of $m=-1,0$ because they give the trivial solution $a_m=0$.

Therefore, for $a_m=x_m(3x_m+2+ t_m)/2$, there is no integer value of $u_m$ that makes $G_{a_m+1,4}(u_m)+2$  zero.

Analogously, for  $a_m=x_m(3x_m+2- t_m)/2$, the absolute values of the real roots of $H(r_m)$ for each value of $s_m$ lie between 0.12 and 19. For $m=-2$, observe that $r_m<0.1$. Consequently, there is no integer value of $u_m$ that makes $G_{a_m+1,4}(u_m)+2$  zero. \EndProof

{\bf Claim 3.} There is no regular graph of odd degree $d\ge3$, girth 9 and cyclic excess.

{\bf Proof of Claim 3.} In this case we proceed as in Claim 2, then  $$H(r_m):=r^3_m(G_{a_m+1,4}(u_m)-2)=r^3_m(G_{a_m+1,4}(x_m+s_m)-2)=0,$$ and $H(r_m)$ is a polynomial in $r_m$ of degree 6.

 For  $a_m=x_m(3x_m+2+ t_m)/2$, the absolute values of the real roots of $H(r_m)$ for each value of $s_m$ lie between 0.05 and 8. For $m=-2$, observe that $a_m=3$, a contradiction. Consequently, there is no integer value of $u_m$ that makes $G_{a_m+1,4}(u_m)-2$  zero.

  Finally, for  $a_m=x_m(3x_m+2- t_m)/2$, the absolute values of the real roots of $H(r_m)$ for each value of $s_m$ lie between 0.12 and 17. Consequently, there is no integer value of $u_m$ that makes $G_{a_m+1,4}(u_m)-2$  zero. \EndProof

  The theorem follows from Claims 2 and 3.
\end{proof}

\subsection{Further non-existence results}
\label{subsec:g5}

\begin{theorem}\label{theo:D3}
There is no regular graph of  degree $d\ge3$, diameter 3 and cyclic defect.
\end{theorem}
\begin{proof}
From Corollary \ref{congruences} it follows that $d-1\equiv0\pmod3$ and that $n\equiv0\pmod3$. In this case we see that $-1$ is an eigenvalue of $B$ with multiplicity 2. Thus,
$G_{d,3}(x)+1=x^3+x^2-(d-1)(2x +1)+1$ must have factors of degree
at most 2,
and therefore an integer root $\lambda$ congruent to 1 modulo 3.
Since $d>1$, we see that $2\lambda+1$ divides $\lambda^3+\lambda ^2+1$, and thus,
$2\lambda+1$ divides 9 (because
$8(\lambda^3+\lambda ^2+1)-9=(2\lambda+1)(4\lambda ^2 +2\lambda -1)$)
and $\lambda\in \{-5,-2,-1,0,1,4\}$. However, from these values only $\lambda=4$ is congruent to 1 modulo 3.

For $\lambda=4$ and $D=3$, we have that $d=10$ and $n=909$. By Proposition \ref{prop:irr} the polynomial$f_9(G_{10,3}(x))$ must be reducible over $\Q[x]$
(see also Observation \ref{elemDiv1}). However,
$$f_9(G_{10,3}(x))=x^9+3x^8-21x^7-74x^6+114x^5+597x^4+160x^3-1488x^2-1920x-701$$
 from where we obtain that $f_9(G_{10,3}(x))$ is irreducible over $\Q[x]$.
\end{proof}

\begin{theorem}\label{theo:OddDiam6}There is no regular graph of odd degree $d\ge3$, diameter $D\equiv0\pmod{6}$, and cyclic defect.
\end{theorem}
\begin{proof}
Since $6|D$, by Corollary \ref{Ddivn}, the order $n$ of these graphs is a multiple of $D$, implying that $n$ is a multiple of 3 and 4. In this case, by Proposition \ref{prop:IntegerRootDef} the polynomial
$G_{d,D}(x)+2$ should have an integer root $\lambda$. On
the other hand, from $6|D$ it follows that $d\equiv1\pmod 6$.
Set $$G_{d,D}(x)+2=x^D+x^{D-1}+(d-1)q(x)+2,$$ where $q(x)$ is a polynomial of degree $D-2$. Thus, $\lambda ^D +\lambda ^{D-1}+2$ should be congruent to 0 modulo 6. But no integer $\lambda$ satisfies $3|(\lambda ^D +\lambda ^{D-1}+2)$.
\end{proof}

\begin{theorem}\label{theo:oddG5}
There is no regular graph of odd degree $d\ge3$, girth 5 and cyclic excess.
\end{theorem}

\begin{proof}
In this case Equation~(\ref{eq:ExM}) takes the form $A^2+A-(d-1)I_n=J_n-B$.
If $d$ is odd then $n=d^2+3\equiv0\pmod{4}$. By Proposition \ref{prop:IntegerRootEx},
there is a simple integer eigenvalue $\lambda$ of $A$ satisfying
\begin{equation}\label{eq1}
\lambda^2 + \lambda -(d-1) =2.
 \end{equation}

As $4|n$, 0 is an eigenvalue of $B$ with multiplicity 2. Therefore, the eigenvalues of $A$
paired with 0 satisfy the equation
\begin{equation}\label{eq2}
x^2 + x -(d-1) =0.
 \end{equation}
Denote by $\lambda_1$ and $\lambda_2$ the roots of Equation~(\ref{eq2}).
If both are eigenvalues of the restriction of $A$ on $\ker(B)$,
the trace is $-1$ (see Corollary \ref{cor:TrMEx}). Therefore, only one of them can be an eigenvalue, say $\lambda_1$, implying that $\lambda_1$ has multiplicity 2 and is an integer.

The discriminant of Equation~(\ref{eq2}) is $4d-3$ and, like the
discriminant $4d+5$ of Equation~(\ref{eq1}), must be a perfect square.
The only pair of perfect squares differing by 8 is \{1,9\}, implying $d=1$,
contradicting the hypothesis $d\ge 3$.
\end{proof}

\begin{theorem}\label{theo:g7}
There is no regular graph of degree $d\equiv0,2\pmod{3}$, girth 7, and cyclic excess.
\end{theorem}
\begin{proof}
Such a graph has  an order multiple of 3. Therefore, the polynomial $G_{d,3}(x)-1$ has a factor of degree 1 or 2, and thus, an integer root $\lambda$. Since $d$ is an integer,  $2\lambda+1$
divides $\lambda ^3+\lambda ^2-1$, and thus, divides 7 (because $8(\lambda^3+\lambda ^2-1)+7=(2\lambda+1)(4\lambda ^2 +2\lambda -1)$). The possible values for $\lambda$  are
$-4$, $-1$, 0 and  3, and the corresponding values for $d$ are 8, 2, 0 and  6. The orders for the
interesting degrees 6 and 8 are $189=3^3\cdot 7$ and $459=3^3\cdot 17$, respectively. But in both cases,
substituting $y=-G_{d,3}(x)$ in $f_9(y)=y^3-3y+1$, we obtain an irreducible polynomial $F^*_{9,d,3}(x)$ of degree 9, contradicting Proposition \ref{prop:irrbis}. Thus, none of these graphs exists.
\end{proof}

\subsection{Computational explorations of graphs of small odd degree with cyclic defect or excess}

In this section we show how to use Corollaries \ref{cor:EvenN} and \ref{cor:EvenNEx}, and the software
Maple\texttrademark\cite{Maple} in order to prove the non-existence of graphs of small degree with cyclic defect or excess. Specifically, we analyze the existence of an integer root in the polynomials
$G_{d,k}(x)\pm2$ for $3\le k\le 20000$ and small degrees. Cubic graphs with cyclic defect or excess are considered in Subsection \ref{sec:Cubic}, while the case of $g=5$ for all graphs of odd degree and cyclic excess
was dealt in Subsection \ref{subsec:g5}. In this subsection we assume $d\ge5$ and $g\ge7$.
\begin{theorem}\label{theo:5Def} For $3\le D\le 20000$
there is no graph of degree 5, diameter $D$ and cyclic defect.
Furthermore, for $7\le g\le 40001$, $g$ odd, there is no graph of degree 5,
girth $g$ and cyclic excess.
\end{theorem}

\begin{proof}
For $3\le k\le 20000$ we analyze the polynomial $G_{5,k}(x)=\pm2$, for $x\in\Z$ and $-4\le x\le 4$. For $x=-4,-2,-1,0,2,3,4$,
we have that $G_{5,k}(x)\equiv0\pmod 4$ if $k\ge3$.

For $3\le k\le 20000$ $G_{5,k}(-3)$ only takes the values $\pm2$ for $k=3,7$,
and in these cases $G_{5,k}(-3)=2$. However, as the order of such graphs
is a multiple of 4. By Corollary \ref{Cor:3,4,6}, both $G_{5,3}(x)$ and $G_{5,7}(x)$
must be reducible over $\Q[x]$, but they are not.

For $3\le k\le 20000$ $G_{5,k}(1)$ takes the values $\pm2$ only for $k=4$, and
then $G_{5,4}(1)=-2$. But, $n\equiv0\pmod{4}$ and $G_{5,4}(x)$ is irreducible over $\Q[x]$, contradicting Corollary \ref{Cor:3,4,6}.
\end{proof}

\begin{theorem}\label{theo:7Def} For any $D\ge3$
there is no graph of degree 7, diameter $D$ and cyclic defect.
Furthermore, for any $g\ge7$, $g$ odd, there is no graph of degree 7,
girth $g$ and cyclic excess.
\end{theorem}

\begin{proof}
 Since $2\sqrt{6}<5$, it suffices to look at $G_{7,k}(x)$ for $x\in\Z$ and $-4\le x\le 4$.
Indeed, $G_{7,3}(x)\ne\pm2$ for $-4\le x\le 4$; for any $k\ge4$ and $x=-4,-3,-2,0,1,2, 4$,
we have that $G_{7,k}(x)\equiv0\pmod{4}$; and  for $k\ge3$ and $x=-1, 3$, it follows
that $G_{7,k}(x)\equiv0\pmod{6}$.
\end{proof}

\begin{theorem}\label{theo:9Def} For $3\le D\le 20000$
there is no graph of degree 9, diameter $D$ and cyclic defect.
Furthermore, for $7\le g\le 40001$, $g$ odd,  there is no graph of degree 9,
girth $g$ and cyclic excess.
\end{theorem}

\begin{proof}
Since $2\sqrt{8}<6$,
it suffices to look at the values of $G_{9,k}(x)$ for $x\in\Z$ and
$-5\le x\le 5$. For $k\ge 3$ and $x\in\{-5,-4,-2,-1,0,2,3, 4\}$,  the value $G_{9,k}(x)$ is a  multiple of 4.
For $3\le k\le 20000$ and $x=-3,1,5$, we have that $G_{9,k}(x)$
never takes the values $\pm2$.
\end{proof}

\begin{theorem}\label{theo:11Def} For $3\le D\le 20000$
there is no graph of degree 11, diameter $D$ and cyclic defect.
Furthermore, for $7\le g\le 40001$, $g$ odd,
 there is no graph of degree 11, girth $g$ and cyclic excess.
\end{theorem}

\begin{proof}
Since $2\sqrt{10}<7$,  it suffices to look at $G_{11,k}(x)$
for  $x\in\Z$ and $-6\le x\le 6$. First, for $k=3$,  $G_{11,3}(x)$ does not take the values 2 or $-2$.
Then, for $k\ge4$ and
$x \in\{ -6,-4,-3,-2,0,\allowbreak 1,2,4,5,6\}$, we have that $G_{11,k}(x)$ is a multiple of 4, while for
$k\ge3$ and $x = -5,-1$, $G_{11,k}(x)$ is a multiple of 10.
Finally, for $4\le k\le 20000$ $G_{11,k}(3)$ never takes the values $\pm2$.
\end{proof}

\begin{theorem}\label{theo:13Def} For $3\le D\le 20000$
there is no graph of degree 13, diameter $D$ and cyclic defect.
Furthermore, for $7\le g\le 40001$, $g$ odd,  there is no graph of degree 13,
girth $g$ and cyclic excess.
\end{theorem}
\begin{proof}
Since $2\sqrt{10}<7$,  it suffices to look at the values of $G_{13,k}(x)$ for  $x\in\Z$
and $-6\le x\le 6$. For $k\ge4$ and $x=-6,-5,-4,-2,-1,0,2,3,4,6$,
we have that $G_{13,k}(x)$ is a multiple of 4, while for
$k\ge3$ and $x=-3, 5$, the polynomial  $G_{13,k}(x)$ is a multiple of 6. For
$3\le k\le20000$ and $x=1$,  $G_{13,k}(x)$ does not take
 the value $2$ or $-2$. Finally, the polynomial $G_{13,3}(x)$
never takes the values $\pm2$ for $-6\le x\le 6$.
\end{proof}

This approach is likely to work for graphs of higher degrees and larger diameters or girths, but its application quickly becomes monotonous and uninteresting.

However,  the aforementioned non-existence results of graphs of odd degree with cyclic defect or excess motivated us to unveil a deeper phenomenon, namely, the finiteness of  such graphs (see Subsection \ref{sec:Finiteness}).

\subsection{Finiteness of graphs of odd degree with cyclic defect or excess}
\label{sec:Finiteness}

In this section we prove the most important results of the paper, namely, the finiteness of all graphs of odd degree $d\ge5$ and cyclic defect or excess (see Theorem \ref{theo:finiteness}),  and the non-existence of cubic graphs with  cyclic defect or excess (see Theorem \ref{theo:3Def}).

The idea behind the proof of Theorem \ref{theo:finiteness} is the following. For any odd degree $d\ge5$ graphs of diameter $k$ and cyclic defect, or graphs of girth $2k+1$ and cyclic excess have an order multiple of 4, implying that the polynomial $G_{d,k}(x)$ must have an algebraic integer of degree at most 2 as a root.  Making use of Equation (\ref{eq:Gcos}) and the fact that any eigenvalue $\lambda(\ne d)$ has the form $2\sqrt{d-1}\cos\alpha$ (with $0<\alpha<\pi$), we show that $\cos\alpha$ must be an algebraic integer of degree at most 4. We then note that if, for a given $d$ and an eigenvalue $\lambda$, Equation (\ref{eq:Gcos}) has at least two values of $k$, then $\alpha$ must be rational. In the case of $\alpha$ being rational and $\cos\alpha$ being an algebraic integer of degree at most 4, we verify that, for $d\ge5$ and all the possible values of $\cos\alpha$,  the polynomial $G_{d,k}(2\sqrt{d-1}\cos\alpha)$ has no algebraic integer of degree at most 2 as a root. This last result implies that for a given $d$ the number of different eigenvalues $\lambda$ of Equation (\ref{eq:Gcos}) represents an upper bound for the number of graphs of degree $d$ and cyclic defect or excess. Finally, we proceed to provide an asymptotic bound for the number of such eigenvalues, knowing that they are algebraic integers of degree 2 lying between $-2\sqrt{d-1}$ and $2\sqrt{d-1}$.

\begin{theorem}\label{theo:finiteness} There are finitely many graphs of odd degree $d\ge5$ and cyclic defect or excess. Furthermore, an asymptotic bound for the number of such graphs is given by $O(\frac{64}{3}d^{3/2})$.
\end{theorem}
\begin{proof}
For graphs of diameter $D=k$ and cyclic defect, and graphs of girth $g=2k+1$ and cyclic excess, if its degree $d$ is odd then its order $n$ is a multiple of 4, which implies, by  Corollary \ref{Cor:3,4,6} (for cyclic defect) and Corollary \ref{Cor:3,4,6bis} (for cyclic excess), that the polynomial $G_{d,k}(x)$  must be reducible over $\Q[x]$.

From Propositions \ref{prop:EigBounds} and \ref{prop:EigBoundsEx} it follows that an eigenvalue $\lambda$  $(\ne d)$ of such graphs has the form $2q\cos\alpha$ with $0<\alpha<\pi$ and $q:=\sqrt{d-1}$, and that $\lambda$ lies between $-2q$ and $2q$. In this case,  because of Equation (\ref{eq:Gcos}) the equation $G_{d,k}(2q\cos\alpha)=0$ implies that
\begin{equation}\label{eq:qsin}q\sin((k+1)\alpha)+\sin(k\alpha)=0.\end{equation}
Also, by Observation \ref{elemDiv1}, such an eigenvalue is an algebraic integer of degree at most 2.

We first claim the following.

{\bf Claim 1.} For a given eigenvalue $\lambda=2q\cos\alpha$, the number $\cos\alpha$ is an algebraic integer of degree at most 4.

{\bf Proof of Claim 1.} Because of Equation (\ref{eq:GdD_UD}), we can expressed Equation (\ref{eq:qsin}) as $$\sin\alpha\left(qU_k(\cos\alpha)+U_{k-1}(\cos\alpha)\right)=0.$$ Then, as $\sin\alpha\ne0$, it follows that \begin{equation}\label{eq:U_k}qU_k(\cos\alpha)+U_{k-1}(\cos\alpha)=0\end{equation}

From Equation (\ref{eq:U_k})  it follows that $(d-1)U^2_k(\cos\alpha)-U^2_{k-1}(\cos\alpha)$ is a polynomial  of degree $2k$ with integer coefficients, having $\cos\alpha$ as a root. Therefore,  $\cos\alpha$ is an algebraic integer of degree at most $2k$.

To see that $\cos\alpha$ is in fact an algebraic integer of degree at most $4$, we need the following facts from algebraic number theory (i) if $\mu$ is an algebraic number of degree $\rho$ then $1/\mu$ is also an algebraic number of degree $\rho$, and (ii) if $\mu$ and $\upsilon$ are algebraic numbers of degree $\rho$ and $\varrho$, respectively, then $\mu\upsilon$ is an algebraic number whose degree divides $\rho\varrho$.

As $\lambda=2q\cos\alpha$ is an algebraic number of degree at most 2 and $q$ is  an algebraic number of degree 2, by the previous facts, $\cos\alpha$ is an algebraic number of degree at most 4.\EndProof

{\bf Claim 2.}  For a given odd  degree $d\ge5$ and an eigenvalue $\lambda$, there is only one value of $k$ satisfying Equation (\ref{eq:qsin}).

{\bf Proof of Claim 2.} We proceed by contradiction, assuming that for a given odd  degree $d\ge5$ and an eigenvalue $\lambda$ ($-2q<\lambda<2q$), there  are at least two values $k_1$ and $k_2$ for which Equation (\ref{eq:qsin}) holds. Observe that in this case $\alpha=\pi r/s$, where $r,s\in \N$. Indeed, assuming that $\sin k\alpha\ne 0$ (for otherwise $\alpha=\pi/2 +p\pi$ with $p\in \N$), Equation (\ref{eq:qsin}) is equivalent to $\cot k\alpha=(-1/q-\cos\alpha)/\sin\alpha$ (since $\sin(k+1)\alpha=\sin k\alpha\cos\alpha+\sin\alpha\cos k\alpha$). If there are two values $k_1$ and $k_2$ for which Equation (\ref{eq:qsin}) holds, then  $$\cot k_1\alpha=\cot k_2\alpha=-\frac{1/q+\cos\alpha}{\sin\alpha}$$
Then, as $\cot x$ is a function with period $\pi$, we have that $(k_1-k_2)\alpha=\pi p$, where $p\in \Z$. In other words, $\alpha/\pi$ is rational.

Therefore, for $0<\alpha<\pi$ we have three cases according to the degree of $\cos\alpha$; see \cite{Jah}.
\begin{itemize}
\item[(i)] If $2\cos\alpha$ is an algebraic integer of degree 1, then $\cos\alpha\in \{-1/2,0,1/2\}$.

\item[(ii)] If $2\cos\alpha$ is an algebraic integer of degree 2, then $$\cos\alpha\in \{-\sqrt3/2,-(1+\sqrt5)/4,-\sqrt2/2,(1-\sqrt5)/4,(-1+\sqrt5)/4,\sqrt2/2\}\cup$$ $$\cup\{(1+\sqrt5)/4,\sqrt3/2\}$$ 

\item[(iii)] If $2\cos\alpha$ is an algebraic integer of degree 4, then  $\alpha=2\pi r/s$ with $r\in \mathbb{N}$ and $s\in \{15,16,20,24,30\}$, or equivalently,
$$\alpha\in\{\pi/15,\pi/12,\pi/10,\pi/8,2\pi/15,4\pi/15,3\pi/10,3\pi/8,5\pi/12,7\pi/15,8\pi/15\}\cup$$ $$\cup\{7\pi/12,5\pi/8,7\pi/10,11\pi/15,13\pi/15,7\pi/8,9\pi/10,11\pi/12,14\pi/15\}$$\end{itemize}
We analyze each case in order, that is, for the aforementioned values of $\cos\alpha$ we look for the values of odd $d$ and $k$ satisfying Equation (\ref{eq:U_k}), or equivalently, $$-q=\frac{U_{k-1}(\cos\alpha)}{U_{k}(\cos\alpha)}.$$
To do this task we sometimes rely on the software Maple\texttrademark\cite{Maple}.

{\bf Case (i)} $\cos\alpha\in \{-1/2,0,1/2\}$.

For $\cos\alpha=-1/2$ and $k\equiv1\pmod3$, we have that $-q=-1$, and thus, $d=2$; for $\cos\alpha=0$ and $k\equiv0\pmod2$, we have $-q=0$; and for $\cos\alpha=1/2$ and $k\equiv1\pmod3$, we have $-q=1$. Therefore, there are no feasible values for $d$ and $k$.

{\bf Case (ii)} $\cos\alpha\in \{-\sqrt3/2,-(1+\sqrt5)/4,-\sqrt2/2,(1-\sqrt5)/4,(-1+\sqrt5)/4,\sqrt2/2,(1+\sqrt5)/4,\sqrt3/2\}$.

The only viable value of $-q$ is $-\sqrt2$, which implies that $d=3$. This case occurs when $\cos\alpha=-\sqrt2/2$ and $k\equiv2\pmod4$.

{\bf Case (iii)} $\alpha=2\pi r/s$ with $r\in \mathbb{N}$ and $s\in \{15,16,20,24,30\}$.

 In this case it can be verified that the only feasible values of $\alpha$ are $5\pi/12$ ($\cos5\pi/12=(-1+\sqrt3)/(2\sqrt2)$) and $11\pi/12$ ($\cos11\pi/12=(-1-\sqrt3)/(2\sqrt2)$). For these values of $\cos\alpha$, we have that $d=3$ ($-q=-\sqrt2$) when $k\equiv9\pmod{12}$.

As a result, when $\alpha/\pi$ is rational there is no odd degree $d\ge5$ satisfying Equation (\ref{eq:qsin}), and thus, the claim follows.\EndProof

Claim 2 also tells us that, for a given odd degree $d$, the number of distinct eigenvalues $\lambda$ is an upper bound for the number of graphs of degree $d$ and cyclic defect and excess.

Recall that, since characteristic polynomials have integer coefficients, if  $\lambda$ is an eigenvalue, so is its conjugate $\lambda^*$.

{\bf Claim 3.} Let $\lambda$ be an eigenvalue of a graph of odd degree and cyclic defect or excess, such that $|\lambda|<2q$ and $|\lambda^*|<2q$. Then, the number of such eigenvalues lying in $(-2q,2q)$ is $O(\frac{64}{3}d^{3/2})$.

{\bf Proof of Claim 3.} We first state a very well known fact about the ring $R_r$ of integers of $\Q(\sqrt r)$, where $r$ is a square-free integer (see \cite[Theorem 1 on pp. 35]{Sam70}): if $r\equiv2,3\pmod4$ then $R_r=\{a+b\sqrt{r}|a,b\in \Z\}$, while if $r\equiv1\pmod4$ then $R_r=\{(u+v\sqrt{r})/2|u,v\in \Z\}$ with $u\equiv v\pmod2$.

By virtue of the previous fact and as $|\lambda|<2q$ and $|\lambda^*|<2q$, we can assume that $\lambda$ has the form either $a+\sqrt{b}$ with $a,b\in\N$ or $(a+\sqrt{b})/2$ with $a,b\in\N$, $a\equiv1\pmod2$ and  $b\equiv1\pmod4$.  In the former case it follows that $0\le a\le2q$ and $1\le b\le(2q-a)^2$, while in the latter we have that $0\le a\le4q$ and $1\le b\le(4q-a)^2$.

An asymptotic bound can be obtained from the remark that the number of pairs $(a,b)$ such that $0\le a\le s$ and that $1\le b\le(s-a)^2$ is at most $(s+1)^3/3$. Indeed,
\begin{align*}
s^2+(s-1)^2+\ldots+(s-\lfloor s\rfloor) &\le \lceil s\rceil^2+(\lceil s\rceil-1)^2+\ldots+(\lceil s\rceil-\lfloor s\rfloor)^2= \\
&=\frac{\lceil s\rceil(\lceil s\rceil+1)(2\lceil s\rceil+1)}{6}<\\
&<\frac{(s+1)^3}{3}.
\end{align*}
The lemma follows from considering the bounds for $a$ and $b$.  \EndProof

Note that the representation of a number $\lambda$ is not unique; for instance, the numbers $1+\sqrt9$ and $2+\sqrt4$ represent the same $\lambda$. However, this detail only makes our bound rougher.

A more careful counting shows that $\frac{64}{3}d^{3/2}$ is indeed an upper bound for the number of algebraic integers of degree 2 lying in $(-2\sqrt{d-1},2\sqrt{d-1})$. As a way of illustration, see that for $d=3$ there are $38<64\sqrt3$ algebraic integers between $-2\sqrt2$ and $2\sqrt2$ while for $d=5$ there are $112<512/3$ algebraic integers lying in $(-4,4)$.

The proof of Theorem \ref{theo:finiteness} follows immediately from Claims 1, 2 and 3.\end{proof}

\subsubsection{Finiteness of cubic graphs with cyclic defect or excess}
\label{sec:Cubic}

Theorem \ref{theo:finiteness} did not settle the finiteness of cubic graphs with cyclic defect or excess. In this subsection we take care of this case.

The proof of Theorem \ref{theo:3Def} first exploits the fact  that for a cubic graph of diameter $k$ and cyclic defect, or a cubic graph of girth $2k+1$ and cyclic excess, the polynomial $G_{3,k}(x)\pm2$ must have an integer root $\rho$ between $-2\sqrt{2}$ and $2\sqrt{2}$. In this direction we prove that $\rho=-1$, and in this case, $G_{3,k}(-1)=2$ and $k\equiv2\pmod{4}$, thus ruling out the existence of graphs with cyclic defect. As the order $n$ of cubic graph with cyclic excess is a multiple of 4,  $G_{3,k}(x)$ has an algebraic root of degree at most 2. In addition, if $k\equiv2\pmod{4}$ then $n\equiv0\pmod3$, so $G_{3,k}(x)\pm1$ must also have algebraic roots of degree 1 or 2. Recall that all these algebraic integers lie on the interval $(-2\sqrt2,2\sqrt2)$. The next step of the proof is to settle that $x=-2$ is the only algebraic integer of degree at most 2 for which the polynomial $G_{3,k}(x)=0$, where $k\equiv 2 \pmod 4$. The conditions that $k\equiv2\pmod4$,  $G_{3,k}(-2)=0$ and $G_{3,k}(-1)=2$ greatly narrow down the numbers $x$ that could make $G_{3,k}(x)=1$; we then prove that the only such numbers are $x=(-1\pm\sqrt{5})/2$ and $x=(-3\pm\sqrt{5})/2$. The proof ends when we prove that, in fact under all the previous conditions, for $x=(-1\pm\sqrt{5})/2$ $x=(-3\pm\sqrt{5})/2$ the polynomial $G_{3,k}(x)\ne1$.

\begin{theorem}\label{theo:3Def} For $D\ge 3$
there is no graph of degree 3, diameter $D$ and cyclic defect.
Furthermore,  for odd $g\ge5$  there is no graph of degree 3,
girth $g$ and cyclic excess.
\end{theorem}

\begin{proof} Suppose, for the sake of contradiction, that there is at least a cubic graph of diameter $D=k\ge3$ and cyclic defect, and at least a graph of girth $g=2k+1\ge5$ and cyclic excess.

Relying on Corollaries \ref{cor:EvenN} and \ref{cor:EvenNEx}, next we discard the values of
 $k\ge3$ and $x\in\{-2,-1,0,1,2\}$ for which $G_{3,k}(x)$ is different from $2$ or $-2$.

Note that $G_{3,k}(2)$, $G_{3,k}(0)$, $G_{3,k}(-2)$ are multiples of 4
if $k\ge4$ and that $G_{3,k}(1)\equiv0\pmod{4}$ if $k\ge2$; this can be checked easily by induction.

For $k=3$ the polynomial $G_{3,3}(x)$ has no factor of degree 1 or 2 to be used with
the eigenvalue 0 of the matrix $B=C_{24}$ (see Corollary \ref{Cor:3,4,6}).
Therefore, there are no cubic graphs of diameter 3 (for cyclic defect) or girth 7 (for cyclic excess).

For $k\ge4$ we can check by induction that $G_{3,k}(-1)\equiv 2\pmod{16}$ if $k$ is even and that
$G_{3,k}(-1)\equiv 10\pmod{16}$
if $k$ is odd. Therefore, from now on we can assume $G_{3,k}(-1)=2$ and $k\equiv0\pmod2$. As a consequence,  there is no cubic graph of diameter $k\ge4$ and cyclic defect.

We now concentrate on cubic graphs of girth $2k+1\ge9$ and cyclic excess, for $k\equiv0\pmod2$.

Computing modulo 32 we see that the  value  $G_{3,k}(-1)=2$  can be attained only if
$k\equiv2 \pmod 4$. In this case,
these graphs have an order multiple of 4 and $3$, and therefore,
we must add the conditions that $G_{3,k}(x)$ and $G_{3,k}(x)\pm1$ have algebraic roots
of degree 1 or 2.

Henceforth, together with the integers $-2,0,1$ and 2, we analyze the set of the 38
algebraic integers of degree 2 between $-2\sqrt2$ and $2\sqrt2$. These numbers are as follows: $\pm\sqrt{u}$ for $u\in\{2,3,5,6,7\}$, $\pm1\pm \sqrt{u}$ for $u\in\{2,3\}$,
$(\pm1\pm \sqrt{u})/2$ for $u\in \{5,13,17,21\}$
and at last $(\pm3\pm\sqrt5)/2$.

We now prove the following.

{\bf Claim 1.} Among all the algebraic integers of degree at most 2 in the interval $(-2\sqrt2,2\sqrt2)$, the polynomial $G_{3,k}(x)$ with $k\equiv 2 \pmod 4$ takes 0 only for $x=-2$.

{\bf Proof of Claim 1.}  The function $G_{3,k}(-2)$ in $k$ is null for $k\equiv 2 \pmod 4$.

The polynomial $G_{3,k}(x)$ never takes the value 0 for $x=0,1,2$, because  $G_{3,k}(0)=-2^{k/2}$, $|G_{3,k}(2)|=2^{k/2+1}$ and $G_{3,k}(1)\equiv 4\pmod8$. These assertions can be proved by induction.

The value $G_{3,k}(x)=0$ is also obtained several times for $x=-1\pm \sqrt3$,
but only if $k\equiv 9\pmod{12}$.

The function $G_{3,k}(\pm\sqrt2)/2^{k/2}$ in $k$ is periodic and never null for $k\equiv 2 \pmod 4$;
the same happens for $1\pm\sqrt3$ and $\pm\sqrt6$.

Consider the algebraic integers $x:=a+b\sqrt{r}$ of degree 2 and odd {\it norm} $N(x)$\footnote{The norm $N$ of quadratic number $x$ is $N(x)=xx^*$.} in the {\it algebraic extension}\footnote{A field extension $L/K$ is called algebraic if every element of $L$ is algebraic over $K$.} $\Q(\sqrt{r})/\Q$. Then, if $x+1$ does not belong to
the {\it principal ideal}\footnote{An {\it ideal} is a subset $I$ of the ring  $R_r$ of algebraic integers  that forms an additive group and has the property that if $\beta\in R_r$ and $\alpha\in I$ then $\beta\alpha\in I$. The {\it principal ideal generated by $\alpha$} with $\alpha\in R_r$ is defined as $\{\alpha\beta|\beta\in R_r\}$.} generated by 2 in  the ring  $R_r$ of algebraic integers, then the polynomial $G_{3,k}$ never enters into that ideal, and thus, never vanishes. This is the case for the numbers $\pm\sqrt3$, $\pm\sqrt7$, $\pm 1\pm\sqrt 2$, $(\pm1\pm\sqrt5)/2$,$(\pm3\pm\sqrt5)/2$, $(\pm1\pm\sqrt{13})/2$, and $(\pm1\pm\sqrt{21})/2$.

For the numbers $(\pm1+\sqrt{17})/2$  we note that for $k\ge6$ the number $G_{3,k}((-1+\sqrt{17})/2)$ has the form $a+b\sqrt{17}$ with $a\equiv{b}\equiv1\pmod2$, and  that for $k\ge10$ the number $G_{3,k}((-1-\sqrt{17})/2)$ has the form $4(a+b\sqrt{17})$ with $a\equiv{b}\equiv1\pmod2$. Furthermore, $G_{3,6}((-1-\sqrt{17})/2)=8$. Therefore, we have ruled out all the numbers $(\pm1\pm\sqrt{17})/2$.

Note that if $G_{3,k}(x)\ne0$ then $G_{3,k}(x^*)\ne0$.

Finally, observe that for $k\ge6$ the number $G_{3,k}(\sqrt5)$ has the form $a+b\sqrt{5}$ with $a\equiv{b}\equiv1\pmod2$. This leaves the numbers $\pm\sqrt5$ out. This completes the proof of the claim.\EndProof

We finalize the proof of the theorem by showing the following two claims.

{\bf Claim 2.} Provided that $k\ge6$ with $k\equiv2\pmod4$,  $G_{3,k}(-2)=0$ and $G_{3,k}(-1)=2$, the only numbers $x$ that could make $G_{3,k}(x)=1$ are $x=(-1\pm\sqrt{5})/2$ and $x=(-3\pm\sqrt{5})/2$.

{\bf Proof of Claim 2.} Set $x:=a+b\sqrt{r}$, then to have simultaneously $G_{3,k}(-2)=0$, $G_{3,k}(-1)=2$ and $G_{3,k}(x)=1$
we must have that $-G_{3,k}(-2) + G_{3,k}(x)=1$ and that $G_{3,k}(-1) - G_{3,k}(x)=1$. Since the coefficients of  $G_{3,k}(x)$ are integers, the former condition means that $x+2$ divides 1 in the ring of integers of
$\Q(\sqrt{r})$, while the latter condition implies that $x+1$ divides 1 in the aforementioned  ring.

These conditions also imply that the norms of $x+2$ and $x+1$ in the algebraic extension $\Q(\sqrt{r})/Q$
must be $1$ or $-1$. Thus, only the numbers $(-1\pm\sqrt 5)/2$ and $(-3\pm\sqrt 5)/2$ satisfy both conditions.\EndProof

{\bf Claim 3.} For  $k\ge6$ with $k\equiv2\pmod4$, the polynomials $G_{3,k}((-1\pm\sqrt 5)/2)$ and $G_{3,k}((-3\pm\sqrt 5)/2)$ never take the value 1.

{\bf Proof of Claim 3.} We first consider the value $(-3+\sqrt 5)/2$, and claim that for $k=2t+4$ with $t\in \N$, the values  $G_{3,k}(\frac{-3+\sqrt 5}2)$ are never integers.

 Observe that $G_{3,2t+4}(x)=(x^2-4)G_{3,2t+2}(x)-4G_{3,2t}(x)$ and that $G_{3,2}((-3+\sqrt5)/{2})=-\sqrt5$. As $N(((-3+\sqrt5)/{2})^2-4)\equiv1\pmod4$ and $N(-\sqrt5)\equiv-1\pmod4$, by induction\footnote{We implicitly use the fact that for algebraic integers $\alpha,\beta$ $N(\alpha\beta)=N(\alpha)N(\beta)$.}, we obtain that $N(G_{3,2t+4}((-3+\sqrt5)/{2}))\equiv-1\pmod4$, which is not the norm of an integer. Recall that the norm of integers is congruent to 0 or 1 modulo 4. This approach also rules out the value $(-3+\sqrt5)/{2}$.

The values of $G_{3,2t+2}((-1+\sqrt 5)/{2})$ are never 1 for $t\ge0$. Indeed,  computing modulo 4 in the ring of integers of $\Q(\sqrt{5})$\footnote{In the ring of the integers $R_r$ of $\Q(\sqrt{r})$, we say that $\alpha$ divides $\beta$, denoted by $\alpha|\beta$, if $\beta/\alpha\in R_r$, and that $\alpha\equiv\beta\pmod\gamma$ if $\gamma|(\alpha-\beta)$.}, we see that $G_{3,2t+2}((-1+\sqrt 5)/{2})\equiv-1\pmod 4 $ if $t$ is multiple of 3, that $G_{3,2t+2}((-1+\sqrt 5)/{2})\equiv\frac{-3+\sqrt5}2$ if
$t\equiv 1\pmod 3$, and that $G_{3,2t+2}((-1+\sqrt 5)/{2})\equiv \frac{-3-\sqrt5}2$ if
$t\equiv 2\pmod 3$. This approach also shows that $G_{3,2t+2}((-1-\sqrt 5)/{2})$ is never 1 for $t\ge0$. This completes the proof of the claim.\EndProof

Combining Claims 1, 2 and 3 the theorem follows.
\end{proof}

An immediate corollary of Theorem of \ref{theo:3Def} is the finiteness of cubic graphs with cyclic defect or excess (Corollary \ref{cor:CubicFiniteness}), settling, in this way, the finiteness of all graphs of odd degree and  cyclic defect or excess.

\begin{corollary}\label{cor:CubicFiniteness} For $k\ge2$, apart from the M\"{o}bius ladder on 8 vertices, there is no cubic graph of diameter $k$ and cyclic defect nor any cubic graph of girth $2k+1$ and cyclic excess.
\end{corollary}

\section{Concluding remarks}
\label{sec:Conclusion}

Using a number of algebraic approaches, we proved the non-existence of infinitely many graphs with cyclic defect or excess, and the finiteness of graphs of odd degree and cyclic defect or excess. While substantial progress in this direction was made through algebraic approaches,  definitive solutions to Problems \ref{ProbRep} and \ref{ProbEx} still seem to be elusive, mainly due to the complexity of the theoretical problems that emerged during our investigation. For instance, the approach which ruled out the existence of cubic graphs with cyclic defect or excess may work for higher degrees, but the complexity of the analysis also increases considerably.

 The condition of having cyclic defect or excess imposes heavy constraints on the structure of graphs with defect or excess 2, so we firmly believe that the M\"{o}bius ladder on 8 vertices is the only such graph, and accordingly, conjecture it.
\begin{conjecture}
Apart from the M\"{o}bius ladder on 8 vertices, there is no graph with cyclic defect or excess.
\end{conjecture}
Furthermore, we think combinatorial approaches have unexplored potential to deal with Problems \ref{ProbRep} and \ref{ProbEx}, so future research should not underestimate them.

\def\cprime{$'$}
\providecommand{\bysame}{\leavevmode\hbox to3em{\hrulefill}\thinspace}
\providecommand{\MR}{\relax\ifhmode\unskip\space\fi MR }
\providecommand{\MRhref}[2]{%
  \href{http://www.ams.org/mathscinet-getitem?mr=#1}{#2}
}
\providecommand{\href}[2]{#2}

\end{document}